\documentclass[11pt,reqno]{amsart}
\usepackage{graphicx}
\usepackage{verbatim}
\usepackage{textcomp}
\usepackage{amssymb}
\usepackage{cite}
\usepackage{amsmath}
\usepackage{latexsym}
\usepackage{amscd}
\usepackage{amsthm}
\usepackage{mathrsfs}
\usepackage{xypic}
\usepackage{bm}
\usepackage{url}
\usepackage{hyperref}

\newtheorem{thm}{Theorem}[section]
\newtheorem{corr}[thm]{Corollary}
\newtheorem{lem}[thm]{Lemma}

\theoremstyle{definition}

\theoremstyle{remark}
\newtheorem{rem}{Remark}[section]
\numberwithin{equation}{section}
\begin{document}
\title[Upper bounds on the first eigenvalue for the $p$-Laplacian]
{Upper bounds on the first eigenvalue for the $p$-Laplacian }
\author{Guangyue Huang}
\address{College of Mathematics and Information Science, Henan Normal
University, Xinxiang, Henan 453007, People's Republic of China}
\email{hgy@henannu.edu.cn }
\author{Zhi Li}
\address{College of Mathematics and Information Science, Henan Normal
University, Xinxiang, Henan 453007, People's Republic of China}
\email{lizhihnsd@126.com} 

\maketitle

\begin{abstract}
In this paper, we establish gradient estimates for positive
solutions to the following equation with respect to the
$p$-Laplacian
$$\Delta_{p}u=-\lambda |u|^{p-2}u$$ with $p>1$ on a given complete
Riemannian manifold. Consequently, we derive upper bound estimates
of the first nontrivial eigenvalue of the $p$-Laplacian.
\end{abstract}

\section{Introduction}

It is well-known that, by using a comparison theorem for the first
Dirichlet eigenvalue of the Laplacian of a geodesic ball on a
complete Riemannian manifold with Ricci curvature bounded from
below, Cheng \cite{Cheng1975} derived an estimate for the bottom of
$L^2$-spectrum on such manifolds by letting the radius of the ball
to infinity. In particular, if the Ricci curvature bounded from
below by $-(n-1)K^2$ for a constant $K$, Cheng obtained the
following upper bound estimate(see Theorem 4.2 in \cite{Cheng1975})
on the first nontrivial eigenvalue of the Laplacian:
\begin{equation}\label{H-In-1}
\lambda_1(\Delta)\leq\frac{(n-1)^2K^2}{4},
\end{equation}

As we know, using the maximum principle to deal with gradient
estimates on Laplacian equations is a power tool in geometric
analysis. For example, Yau \cite{Yau1975} showed that every positive
or bounded solution to the heat equation is constant. Li
\cite{LiD2005} derived several Liouville type theorems with respect
to the weighted Laplacian. Generalizing the Theorem 4.2 of Cheng
\cite{Cheng1975}, Wang \cite{Wang2010} and Wu \cite{Wu2010} obtained
upper bound estimates of weighted Laplacian with respect to the
Bakry-\'{E}mery Ricci curvature, respectively.

It is a natural question to derive upper bound estimates of the
first eigenvalue of the $p$-Laplacian. The $p$-Laplacian on an
$n$-dimensional Riemannian manifold $(M^n,g)$ is defined by
$$\Delta_{p}u:={\rm div}(|\nabla u|^{p-2}\nabla u),\ \ \ \ {\rm for }\ u\in W^{1,p}(M^n)$$ with
$p>1$, which can be seen as a generalization of the Laplacian. The
first nonzero eigenvalue of $\Delta_{p}$ is defined by the
variational characterization
\begin{equation}\label{H-In-2}\lambda_1(\Delta_{p})=\inf\limits_{u\in W^{1,p}(M^n)}\Bigg\{\int\limits_{M^n}|\nabla u|^p\ \Bigg|\ \
\int\limits_{M^n}|u|^p=1,\,
\int\limits_{M^n}|u|^{p-2}u=0\Bigg\}.\end{equation} It has been
proved by Serrin and Tolksdorf in \cite{Tolksdorf1984,Serrin1964}
that the infimum defined in \eqref{H-In-2} is attained by a
$C^{1,\alpha}$ eigenfunction $u$ which satisfies Euler-Lagrange
equation
\begin{equation}\label{In-1}
\Delta_{p}u=-\lambda |u|^{p-2}u.
\end{equation}
In this paper, we first derive locally gradient estimates of
positive solutions to \eqref{In-1} on a geodesic ball. Then main
results of this paper follow from letting the radius of the geodesic
ball to the infinity. Our main results are as follows:

\begin{thm}\label{ThmInt1}
Let $(M, g)$ be an $n$-dimensional complete Riemannian manifold with
the sectional curvature $K_M\geq-K^2$ for some positive constant
$K$. Let $u$ be a positive solution to \eqref{In-1} with $p>2$, then
for $h=(p-1)\log u$,
\begin{equation}\label{Th-1}\aligned
|\nabla h|^p\leq&\frac{(n-1)(p-1)^2}{p}\Bigg\{-\frac{p}{n-1}(p-1)^{p-1}\lambda+(n-1)^{p-1}K^p\Big(\frac{p-1}{p}\Big)^{p-2}\\
&+\Bigg[\frac{p^3(p-1)^{2(p-1)}}{(n-1)^2(p-2)}\lambda^2-\frac{2K^p(p-1)^{2p-3}(n-1)^{p-2}}{p^{p-3}}\lambda\\
& +K^{2p}(n-1)^{2(p-1)}\Big(\frac{p-1}{p}\Big)^{2(p-2)}
\Bigg]^{\frac{1}{2}}\Bigg\},
\endaligned\end{equation} and the first
eigenvalue $\lambda_1$ of eigenvalue problem \eqref{In-1} satisfies
\begin{equation}\label{Th-2}\lambda_1(\Delta_{p})\leq
\Big(\frac{(n-1)K}{p}\Big)^p.\end{equation}

\end{thm}

\begin{thm}\label{ThmInt2}
Let $(M, g)$ be an $n$-dimensional complete Riemannian manifold with
the sectional curvature $K_M\geq-K^2$ for some positive constant
$K$. Let $u$ be a positive solution to \eqref{In-1} with $1<p<2$,
then for $h=(p-1)\log u$,
\begin{equation}\label{2Th-1}\aligned
|\nabla h|^p\leq&\frac{n-1}{p}
\Bigg\{-\frac{p}{n-1}(p-1)^{p-1}\lambda+(n-1)^{p-1}p^{2-p}(p-1)^{p-1}K^p
\\
&+\Bigg[-\frac{2(p-1)^{2(p-1)}p^{3-p}K^p}{(n-1)^{2-p}}\lambda+2(n-1)^{2(p-1)}(p-1)^{2(p-1)}p^{3-2p}K^{2p}\Bigg]^{\frac{1}{2}}\Bigg\}
\endaligned\end{equation} and the first
eigenvalue $\lambda_1$ of eigenvalue problem \eqref{In-1} satisfies
\begin{equation}\label{2Th-2}\lambda_1(\Delta_{p})\leq
\Big(\frac{(n-1)K}{p}\Big)^p.\end{equation}

\end{thm}

In particular, we obtain the following estimate to positive
$p$-harmonic function $u$ satisfying
\begin{equation}\label{Harmo-In-1}
\Delta_{p}u=0
\end{equation} by letting $\lambda=0$ in \eqref{Th-1} and \eqref{2Th-1}, respectively:

\begin{corr}\label{corr1}
Let $(M, g)$ be an $n$-dimensional complete Riemannian manifold with
the sectional curvature $K_M\geq-K^2$ for some constant $K$. Let $u$
be a positive solution to \eqref{Harmo-In-1} and $h=(p-1)\log u$.
Then,

(1) for $p>2$,
\begin{equation}\label{corr-1}\aligned
|\nabla h|^p\leq&\frac{2}{p^{p-1}}[(n-1)(p-1)K]^p;
\endaligned\end{equation}

(2) for $1<p<2$,
\begin{equation}\label{2corr-1}\aligned
|\nabla
h|^p\leq&\Big(1+\sqrt{\frac{2}{p}}\Big)\Big(\frac{p-1}{p}\Big)^{p-1}[(n-1)K]^p.
\endaligned\end{equation}
\end{corr}

\begin{rem} Letting $p\rightarrow2$ in \eqref{Th-2} or
\eqref{2Th-2}, we obtain the upper bound \eqref{H-In-1} of the first
nontrivial eigenvalue of the Laplacian of Cheng \cite{Cheng1975}. We
also obtain that
\begin{equation}\label{Int11}
|\nabla h|\leq(n-1)K
\end{equation} by letting $p\rightarrow2$ in \eqref{corr-1} or in
\eqref{2corr-1}. In particular, the constant $(n-1)K$ is sharp in
light of the results of \cite{LiYau1986,Yau1975} for $p=2$. As a
direct consequence, we derive from \eqref{corr-1} and
\eqref{2corr-1} that if $(M, g)$ is a complete Riemannian manifold
with nonnegative sectional curvature, then any positive
$p$-harmonic($p>1$) function must be a constant.

\end{rem}

\begin{rem} It has been shown in \cite{WangZhu2012} that the
estimate given in \eqref{2Th-2} is sharp when $1<p<2$. Similarly, we
can prove that the estimate given \eqref{Th-2} is also sharp when
$p>2$. Some related results for the lower bound of the first
nontrivial eigenvalue of $p$-Laplacian can be found in
\cite{zhanghc07,LFW09,Val1,Matei,Wang2016,WangLi2016} and the
references therein.

\end{rem}

%\begin{ack}
%The research of authors is supported by NSFC(No. 11371018,
%11401179), Henan Provincial Backbone Teacher(No. 2013GGJS-057), and
%Henan Provincial Education Department(No. 14B110017).
%\end{ack}

\section{The Proof of the Case $p>2$}

We define
\begin{equation}\label{In-2}
h=(p-1)\log u.
\end{equation} Then from \eqref{In-1}, we obtain that $h$ satisfies
\begin{equation}\label{In-3}
\Delta_{p}h=-|\nabla h|^p-(p-1)^{p-1}\lambda.
\end{equation} Next, we give a key lemma which
has been proved by Wang and Li \cite{WangLi2016}(see Lemma 3.2 in
\cite{WangLi2016}):

\begin{lem}\label{2lemma1} We introduce the elliptic operator
$\mathcal{L}$ defined by
$$\mathcal{L}={\rm div}(|\nabla h|^{p-2}A\nabla\cdot)$$ with
$A=Id+(p-2)\frac{\nabla h\otimes\nabla h}{|\nabla h|^{2}}$.
 Then, for $p>1$ and
 \begin{equation}\label{In-4} G=|\nabla h|^p,
\end{equation} we have
\begin{equation}\label{In-5}\aligned
\mathcal{L}(G)\geq&\frac{p}{n-1}G^2+\frac{2p}{n-1}(p-1)^{p-1}\lambda
G+\frac{p(p-1)^{2(p-1)}}{n-1}\lambda^2\\
&-p(n-1)K^2G^{\frac{2(p-1)}{p}}+\frac{\alpha}{p}\frac{|\nabla
G|^2}{G^{\frac{2}{p}}}
\\
&+\Bigg[\frac{2(p-1)}{n-1}\Big(1+\lambda
G^{-1}(p-1)^{p-1}\Big)-p\Bigg]\langle\nabla h,\nabla G\rangle
|\nabla h|^{p-2},
\endaligned\end{equation} where
$\alpha:=\min\{2(p-1),\frac{n(p-1)^2}{n-1}\}$.
\end{lem}

Let $\theta$ be a cut-off function satisfying $\theta(t)=1$ for
$0\leq t\leq \frac{1}{2}$ and $\theta(t)=0$ for $t\geq1$ such that
$$\frac{(\theta')^2}{\theta}\leq10,\ \ \
\theta''\geq-10\theta\geq-10.$$ Define the function
$\varphi:M^n\rightarrow \mathbb{R}$ by
$$\varphi(x)=\theta\Big(\frac{r(x)}{R}\Big),$$ then a direct
calculation shows
\begin{equation}\label{H-In-6}
\frac{|\nabla\varphi|^2}{\varphi}\leq\frac{10}{R^2},
\end{equation} and
\begin{equation}\label{HH-In-6}\aligned
A^{ij}\varphi_{ij}=&\Delta\varphi+(p-2)\frac{\varphi_{ij}h^ih^j}{|\nabla h|^2}\\
\leq&-80(n+p-2)\frac{1+2KR}{R^2}-40\max\{p-1,1\}\frac{1}{R^2}.
\endaligned\end{equation}  For the proof of \eqref{HH-In-6}, we refer to
\cite{KoNi}.

Let $\tilde{G}=\varphi G$. Next we will apply the maximum principle
to $\tilde{G}$ on $B_p(R)$. We assume that $\tilde{G}=\varphi G$
achieve its maximal value at $x_0$. Then at the point $x_0$, we have
\begin{equation}\label{In-7}
\nabla G=-G\frac{\nabla\varphi}{\varphi}.
\end{equation}
It is easy to see that
\begin{equation}\label{In-6}\aligned
\mathcal{L}(\tilde{G})=&\varphi\mathcal{L}
G+G\mathcal{L}\varphi+2\langle\nabla\varphi,A(\nabla
G)\rangle|\nabla h|^{p-2}.
\endaligned\end{equation}

Next, we estimate each terms on the right hand of \eqref{In-6} for
$p>2$. From \eqref{In-5}, we have
\begin{equation}\label{In-8}\aligned
\varphi\mathcal{L} G \geq&\frac{p}{n-1}\
\frac{\tilde{G}^2}{\varphi}+\frac{2p}{n-1}(p-1)^{p-1}\lambda
\tilde{G}+\frac{p(p-1)^{2(p-1)}}{n-1}\lambda^2\varphi\\
&-p(n-1)K^2\frac{\tilde{G}^{\frac{2(p-1)}{p}}}{\varphi^{\frac{p-2}{p}}}
+\frac{\alpha}{p}\frac{|\nabla\varphi|^2}{\varphi}\Big(\frac{\tilde{G}}{\varphi}\Big)^{\frac{2(p-1)}{p}}
\\
&-\frac{\tilde{G}}{\varphi}\Bigg(\frac{2(p-1)}{n-1}
-p\Bigg)\langle\nabla h,\nabla \varphi\rangle
|\nabla h|^{p-2}\\
&-\frac{2}{n-1}(p-1)^p\lambda \langle\nabla h,\nabla \varphi\rangle
|\nabla h|^{p-2}\\
\geq&\frac{p}{n-1}\
\frac{\tilde{G}^2}{\varphi}+\frac{2p}{n-1}(p-1)^{p-1}\lambda
\tilde{G}+\frac{p(p-1)^{2(p-1)}}{n-1}\lambda^2\varphi\\
&-p(n-1)K^2\frac{\tilde{G}^{\frac{2(p-1)}{p}}}{\varphi^{\frac{p-2}{p}}}
+\frac{\alpha}{p}\frac{|\nabla\varphi|^2}{\varphi}\Big(\frac{\tilde{G}}{\varphi}\Big)^{\frac{2(p-1)}{p}}
\\
&-\Big|\frac{2(p-1)}{n-1}
-p\Big||\nabla\varphi|\Big(\frac{\tilde{G}}{\varphi}\Big)^{\frac{2p-1}{p}}-\frac{2}{n-1}(p-1)^p\lambda
|\nabla\varphi|\Big(\frac{\tilde{G}}{\varphi}\Big)^{\frac{p-1}{p}}.
\endaligned\end{equation} By the definition of $A$, we obtain
\begin{equation}\label{In-9}\aligned
2\langle\nabla\varphi,A(\nabla G)\rangle|\nabla
h|^{p-2}=&2\Bigg[\langle\nabla\varphi, \nabla G
\rangle+(p-2)\frac{\langle\nabla \varphi,\nabla
h\rangle\langle\nabla h,\nabla G\rangle}{|\nabla
h|^{2}}\Bigg]|\nabla h|^{p-2}\\
=&-2G\Bigg[\frac{|\nabla\varphi|^2}{\varphi}+(p-2)\frac{\langle\nabla
h,\nabla \varphi\rangle^2}{\varphi|\nabla h|^{2}}\Bigg]|\nabla
h|^{p-2}\\
\geq&-2G\Bigg[\frac{|\nabla\varphi|^2}{\varphi}+(p-2)\frac{|\nabla\varphi|^2}{\varphi}\Bigg]|\nabla
h|^{p-2}\\
=&-2(p-1)\frac{|\nabla\varphi|^2}{\varphi}\Big(\frac{\tilde{G}}{\varphi}\Big)^{\frac{2(p-1)}{p}}.
\endaligned\end{equation}
On the other hand, due to
$$\langle\nabla\varphi,\nabla |\nabla
h|^{p-2}\rangle=\Big(\frac{2}{p}-1\Big)|\nabla
h|^{p-2}\frac{|\nabla\varphi|^2}{\varphi}$$ and
$$|\nabla
h|^{p-2}\langle\nabla\frac{\langle\nabla\varphi,\nabla
h\rangle}{|\nabla h|^2} ,\nabla h\rangle=|\nabla
h|^{p-2}\Big(\frac{\varphi_{ij}h_ih_j}{|\nabla
h|^2}-\frac{1}{p}\frac{|\nabla\varphi|^2}{\varphi}
+\frac{2}{p}\frac{\langle\nabla\varphi,\nabla h\rangle^2}{|\nabla
h|^2\varphi}\Big),$$ we have
\begin{equation}\label{In-10}\aligned
\mathcal{L}\varphi=&{\rm div}(|\nabla
h|^{p-2}A(\nabla\varphi))\\
=&|\nabla
h|^{p-2}\Delta\varphi+(p-2)\frac{\langle\nabla\varphi,\nabla
h\rangle}{|\nabla h|^2}\Delta_ph+\langle\nabla\varphi,\nabla |\nabla
h|^{p-2}\rangle\\
&+(p-2)|\nabla
h|^{p-2}\langle\nabla\frac{\langle\nabla\varphi,\nabla
h\rangle}{|\nabla h|^2} ,\nabla h\rangle\\
=&\Big(\frac{\tilde{G}}{\varphi}\Big)^{\frac{p-2}{p}}\Delta\varphi-(p-2)\frac{\langle\nabla\varphi,\nabla
h\rangle}{|\nabla
h|^2}\Big(G+(p-1)^{p-1}\lambda\Big)\\
&+(p-2)\frac{\varphi_{ij}h_ih_j}{|\nabla
h|^2}\Big(\frac{\tilde{G}}{\varphi}\Big)^{\frac{p-2}{p}}
+2\Big(\frac{2}{p}-1\Big)\frac{|\nabla\varphi|^2}{\varphi}\Big(\frac{\tilde{G}}{\varphi}\Big)^{\frac{p-2}{p}}\\
&+\frac{2(p-2)}{p}\frac{\langle\nabla\varphi,\nabla
h\rangle^2}{|\nabla
h|^2\varphi}\Big(\frac{\tilde{G}}{\varphi}\Big)^{\frac{p-2}{p}}\\
=&\varphi_{ij}A^{ij}\Big(\frac{\tilde{G}}{\varphi}\Big)^{\frac{p-2}{p}}-(p-2)\frac{\langle\nabla\varphi,\nabla
h\rangle}{|\nabla
h|^2}\Big(G+(p-1)^{p-1}\lambda\Big)\\
&+2\Big(\frac{2}{p}-1\Big)\frac{|\nabla\varphi|^2}{\varphi}\Big(\frac{\tilde{G}}{\varphi}\Big)^{\frac{p-2}{p}}
+\frac{2(p-2)}{p}\frac{\langle\nabla\varphi,\nabla
h\rangle^2}{|\nabla
h|^2\varphi}\Big(\frac{\tilde{G}}{\varphi}\Big)^{\frac{p-2}{p}},
\endaligned\end{equation}
which shows that
\begin{equation}\label{In-11}\aligned
G\mathcal{L}\varphi=&\varphi_{ij}A^{ij}\Big(\frac{\tilde{G}}{\varphi}\Big)^{\frac{2(p-1)}{p}}-(p-2)\langle\nabla\varphi,\nabla
h\rangle\Big(\frac{\tilde{G}}{\varphi}\Big)^{\frac{2(p-1)}{p}}\\
&-(p-2)(p-1)^{p-1}\lambda\langle\nabla\varphi,\nabla
h\rangle\Big(\frac{\tilde{G}}{\varphi}\Big)^{\frac{p-2}{p}}\\
&+2\Big(\frac{2}{p}-1\Big)\frac{|\nabla\varphi|^2}{\varphi}\Big(\frac{\tilde{G}}{\varphi}\Big)^{\frac{2(p-1)}{p}}
+\frac{2(p-2)}{p}\frac{\langle\nabla\varphi,\nabla
h\rangle^2}{|\nabla
h|^2\varphi}\Big(\frac{\tilde{G}}{\varphi}\Big)^{\frac{2(p-1)}{p}}\\
\geq&\varphi_{ij}A^{ij}\Big(\frac{\tilde{G}}{\varphi}\Big)^{\frac{2(p-1)}{p}}
-(p-2)|\nabla\varphi|\Big(\frac{\tilde{G}}{\varphi}\Big)^{\frac{2p-1}{p}}\\
&-(p-2)(p-1)^{p-1}\lambda|\nabla\varphi|\Big(\frac{\tilde{G}}{\varphi}\Big)^{\frac{p-1}{p}}
+2\Big(\frac{2}{p}-1\Big)\frac{|\nabla\varphi|^2}{\varphi}\Big(\frac{\tilde{G}}{\varphi}\Big)^{\frac{2(p-1)}{p}}.
\endaligned\end{equation}
Putting \eqref{In-8}, \eqref{In-9} and \eqref{In-11} into
\eqref{In-6} yields
\begin{equation}\label{In-12}\aligned
\mathcal{L}(\tilde{G})\geq&\frac{p}{n-1}\
\frac{\tilde{G}^2}{\varphi}+\frac{2p}{n-1}(p-1)^{p-1}\lambda
\tilde{G}+\frac{p(p-1)^{2(p-1)}}{n-1}\lambda^2\varphi\\
&+\Bigg[-p(n-1)K^2\varphi
+\frac{\alpha}{p}\frac{|\nabla\varphi|^2}{\varphi}-2(p-1)\frac{|\nabla\varphi|^2}{\varphi}\\
&+\varphi_{ij}A^{ij}+2\Big(\frac{2}{p}-1\Big)\frac{|\nabla\varphi|^2}{\varphi}\Bigg]\Big(\frac{\tilde{G}}{\varphi}\Big)^{\frac{2(p-1)}{p}}\\
&-\Bigg[\Big|\frac{2(p-1)}{n-1}
-p\Big|+(p-2)\Bigg]|\nabla\varphi|\Big(\frac{\tilde{G}}{\varphi}\Big)^{\frac{2p-1}{p}}\\
&-\frac{(n+1)p-2n}{n-1}(p-1)^{p-1}\lambda|\nabla\varphi|\Big(\frac{\tilde{G}}{\varphi}\Big)^{\frac{p-1}{p}}.
\endaligned\end{equation}
Using the Cauchy inequality, for any positive constant
$\varepsilon_1$, we have
\begin{equation}\label{In-13}\aligned
-\frac{(n+1)p-2n}{n-1}&(p-1)^{p-1}\lambda|\nabla\varphi|\Big(\frac{\tilde{G}}{\varphi}\Big)^{\frac{p-1}{p}}\\
\geq&-\sqrt{10}\frac{(n+1)p-2n}{n-1}(p-1)^{p-1}\frac{\lambda}{R}\Big(\frac{\tilde{G}}{\varphi}\Big)^{\frac{p-1}{p}}\varphi^{\frac{1}{2}}\\
\geq&-\varepsilon_1\lambda^2
-\frac{5(p-1)^{2(p-1)}\frac{[(n+1)p-2n]^2}{(n-1)^2}}{2\varepsilon_1R^2}\varphi\Big(\frac{\tilde{G}}{\varphi}\Big)^{\frac{2(p-1)}{p}},
\endaligned\end{equation} and for or any positive constant
$\varepsilon_2$,
\begin{equation}\label{In-14}\aligned
-&\Bigg[\Big|\frac{2(p-1)}{n-1}
-p\Big|+(p-2)\Bigg]|\nabla\varphi|\Big(\frac{\tilde{G}}{\varphi}\Big)^{\frac{2p-1}{p}}\\
=&-\Bigg[\Big|\frac{2(p-1)}{n-1}
-p\Big|+(p-2)\Bigg]|\nabla\varphi|\Big(\frac{\tilde{G}}{\varphi}\Big)^{\frac{p-1}{p}}\Big(\frac{\tilde{G}}{\varphi}\Big)\\
\geq&-\varepsilon_2\Big(\frac{\tilde{G}}{\varphi}\Big)^2-\Bigg[\Big|\frac{2(p-1)}{n-1}
-p\Big|+(p-2)\Bigg]^2\frac{5\varphi}{2\varepsilon_2R^2}\Big(\frac{\tilde{G}}{\varphi}\Big)^{\frac{2(p-1)}{p}}.
\endaligned\end{equation}
Inserting \eqref{In-13} and \eqref{In-14} into \eqref{In-12}, we
derive
\begin{equation}\label{In-15}\aligned
\mathcal{L}(\tilde{G})\geq&\Big[\frac{p}{n-1}\varphi-\varepsilon_2\Big]\Big(\frac{\tilde{G}}{\varphi}\Big)^2
+\frac{2p}{n-1}(p-1)^{p-1}\lambda\varphi
\frac{\tilde{G}}{\varphi}+\Big[\frac{p(p-1)^{2(p-1)}}{n-1}\varphi-\varepsilon_1\Big]\lambda^2\\
&+\Bigg\{-p(n-1)K^2\varphi
+\frac{\alpha}{p}\frac{|\nabla\varphi|^2}{\varphi}-2(p-1)\frac{|\nabla\varphi|^2}{\varphi}+\varphi_{ij}A^{ij}\\
&+2\Big(\frac{2}{p}-1\Big)\frac{|\nabla\varphi|^2}{\varphi}-\frac{5(p-1)^{2(p-1)}\frac{[(n+1)p-2n]^2}{(n-1)^2}\varphi}{2\varepsilon_1R^2}\\
&-\Bigg[\Big|\frac{2(p-1)}{n-1}
-p\Big|+(p-2)\Bigg]^2\frac{5\varphi}{2\varepsilon_2R^2}\Bigg\}\Big(\frac{\tilde{G}}{\varphi}\Big)^{\frac{2(p-1)}{p}}\\
=&\Big[\frac{p}{n-1}\varphi-\varepsilon_2\Big]\Big(\frac{\tilde{G}}{\varphi}\Big)^2+\frac{2p}{n-1}(p-1)^{p-1}\lambda\varphi
\frac{\tilde{G}}{\varphi}+\Big[\frac{p(p-1)^{2(p-1)}}{n-1}\varphi-\varepsilon_1\Big]\lambda^2\\
&+\Bigg\{-p(n-1)K^2\varphi
-\frac{2p^2-(\alpha+4)}{p}\frac{|\nabla\varphi|^2}{\varphi}\\
&+\varphi_{ij}A^{ij}-\frac{5(p-1)^{2(p-1)}\frac{[(n+1)p-2n]^2}{(n-1)^2}\varphi}{2\varepsilon_1R^2}\\
&-\Bigg[\Big|\frac{2(p-1)}{n-1}
-p\Big|+(p-2)\Bigg]^2\frac{5\varphi}{2\varepsilon_2R^2}\Bigg\}\Big(\frac{\tilde{G}}{\varphi}\Big)^{\frac{2(p-1)}{p}}\\
\geq&\Big[\frac{p}{n-1}\varphi-\varepsilon_2\Big]\Big(\frac{\tilde{G}}{\varphi}\Big)^2+\frac{2p}{n-1}(p-1)^{p-1}\lambda\varphi
\frac{\tilde{G}}{\varphi}+\Big[\frac{p(p-1)^{2(p-1)}}{n-1}\varphi-\varepsilon_1\Big]\lambda^2\\
&-D\Big(\frac{\tilde{G}}{\varphi}\Big)^{\frac{2(p-1)}{p}},
\endaligned\end{equation}
where
$$\aligned D=&p(n-1)K^2\varphi
+\frac{2p^2-(\alpha+4)}{p}\frac{10}{R^2}\\
&+80(n+p-2)\frac{1+KR}{R^2}+\frac{40(p-1)}{R^2}+\frac{5(p-1)^{2(p-1)}\frac{[(n+1)p-2n]^2}{(n-1)^2}\varphi}{2\varepsilon_1R^2}\\
&+\Bigg[\Big|\frac{2(p-1)}{n-1}
-p\Big|+(p-2)\Bigg]^2\frac{5\varphi}{2\varepsilon_2R^2}.\endaligned$$
Here we used \eqref{H-In-6} and \eqref{HH-In-6} in the last
inequality of \eqref{In-15}. By virtue of Holder inequality again,
one has
\begin{equation}\label{In-16}\aligned
D\Big(\frac{\tilde{G}}{\varphi}\Big)^{\frac{2(p-1)}{p}}
=&\Bigg[\varepsilon_3^{\frac{p-2}{p}}\Big(\frac{\tilde{G}}{\varphi}\Big)^{\frac{2p-4}{p}}\Bigg]
\Bigg[\frac{D}{\varepsilon_3^{\frac{p-2}{p}}}\Big(\frac{\tilde{G}}{\varphi}\Big)^{\frac{2}{p}}\Bigg]\\
\leq&\frac{p-2}{p}\varepsilon_3\Big(\frac{\tilde{G}}{\varphi}\Big)^2+\frac{2}{p}
\frac{D^{\frac{p}{2}}}{\varepsilon_3^{\frac{p-2}{2}}}\frac{\tilde{G}}{\varphi},
\endaligned\end{equation} where $\varepsilon_3$ is a positive constant to be
determined. Thus, \eqref{In-15} can be written as
\begin{equation}\label{In-17}\aligned
\mathcal{L}(\tilde{G})\geq&
\Big[\frac{p}{n-1}\varphi-\varepsilon_2-\frac{p-2}{p}\varepsilon_3\Big]\Big(\frac{\tilde{G}}{\varphi}\Big)^2\\
&+\Big[\frac{2p}{n-1}(p-1)^{p-1}\lambda\varphi-\frac{2}{p}
\frac{D^{\frac{p}{2}}}{\varepsilon_3^{\frac{p-2}{2}}}\Big]
\frac{\tilde{G}}{\varphi}+\Big[\frac{p(p-1)^{2(p-1)}}{n-1}\varphi-\varepsilon_1\Big]\lambda^2,
\endaligned\end{equation}
which gives at the point $x_0$,
\begin{equation}\label{In-18}\aligned
0\geq&
\Big[\frac{p}{n-1}\varphi-\varepsilon_2-\frac{p-2}{p}\varepsilon_3\Big]\tilde{G}^2\\
&+\Big[\frac{2p}{n-1}(p-1)^{p-1}\lambda\varphi^2-\frac{2}{p}
\frac{D^{\frac{p}{2}}}{\varepsilon_3^{\frac{p-2}{2}}}\varphi\Big]
\tilde{G}+\Big[\frac{p(p-1)^{2(p-1)}}{n-1}\varphi^3-\varepsilon_1\varphi^2\Big]\lambda^2.
\endaligned\end{equation}

In order to obtain the bound of $\tilde{G}$ by using the maximum
principle for \eqref{In-18}, it is sufficient to choose the
coefficient of $\tilde{G}^2$ in \eqref{In-18} is positive by taking
$\varepsilon_2,\varepsilon_3$ small enough. Since $x_0$ is the
maximum point of the function $\tilde{G}$ and $\varphi=1$ on
$B_p(\frac{R}{2})$,
$$\varphi(x_0)|\nabla h|^p(x_0)\geq\sup\limits_{B_p(\frac{R}{2})}|\nabla h|^p(x).$$
On the other hand, using the fact that
$$\varphi(x_0)|\nabla h|^p(x_0)\leq\varphi(x_0)\sup\limits_{B_p(R)}|\nabla h|^p(x),$$
it is easy to see that
$$\sigma(R)\leq\varphi(x_0)\leq1,$$ where
$\sigma(R)$ is defined by
$$\sigma(R):=\frac{\sup\limits_{B_p(\frac{R}{2})}|\nabla h|^p(x)}{\sup\limits_{B_p(R)}|\nabla
h|^p(x)}.$$ Thus, at the point $x_0$, \eqref{In-18} gives
\begin{equation}\label{In-19}\aligned
0\geq&
\Big[\frac{p}{n-1}\sigma(R)-\varepsilon_2-\frac{p-2}{p}\varepsilon_3\Big]\tilde{G}^2\\
&+\Big[\frac{2p}{n-1}(p-1)^{p-1}\lambda\sigma(R)^2-\frac{2}{p}
\frac{D^{\frac{p}{2}}}{\varepsilon_3^{\frac{p-2}{2}}}\Big]
\tilde{G}+\Big[\frac{p(p-1)^{2(p-1)}}{n-1}\sigma(R)^3-\varepsilon_1\Big]\lambda^2,
\endaligned\end{equation}
 as long as
\begin{equation}\label{H-In-19}\frac{p}{n-1}\sigma(R)-\varepsilon_2-\frac{p-2}{p}\varepsilon_3>0.\end{equation}
By using the inequality
$$ax^2+bx+c\leq0$$ with $a>0$, then
$$x\leq\frac{-b+\sqrt{b^2-4ac}}{2a},$$ we obtain from
\eqref{In-19}
\begin{equation}\label{HHIn-20}\aligned\Big(&\frac{2p}{n-1}(p-1)^{p-1}\lambda\sigma(R)^2-\frac{2}{p}
\frac{D^{\frac{p}{2}}}{\varepsilon_3^{\frac{p-2}{2}}}\Big)^2\\
&-4\Big(\frac{p}{n-1}\sigma(R)
-\varepsilon_2-\frac{p-2}{p}\varepsilon_3\Big)\Big(\frac{p(p-1)^{2(p-1)}}{n-1}\sigma(R)^3-\varepsilon_1\Big)
\lambda^2\geq0\endaligned\end{equation} and
\begin{equation}\label{In-20}\aligned
\tilde{G}(x_0)\leq&
\frac{1}{2\Big[\frac{p}{n-1}\sigma(R)-\varepsilon_2-\frac{p-2}{p}\varepsilon_3\Big]}\Bigg\{-\Big[\frac{2p}{n-1}(p-1)^{p-1}\lambda\sigma(R)^2-\frac{2}{p}
\frac{D^{\frac{p}{2}}}{\varepsilon_3^{\frac{p-2}{2}}}\Big]\\
&+\Bigg[\Big(\frac{2p}{n-1}(p-1)^{p-1}\lambda\sigma(R)^2-\frac{2}{p}
\frac{D^{\frac{p}{2}}}{\varepsilon_3^{\frac{p-2}{2}}}\Big)^2\\
&-4\Big(\frac{p}{n-1}\sigma(R)
-\varepsilon_2-\frac{p-2}{p}\varepsilon_3\Big)\Big(\frac{p(p-1)^{2(p-1)}}{n-1}\sigma(R)^3-\varepsilon_1\Big)
\lambda^2\Bigg]^{\frac{1}{2}}\Bigg\}.
\endaligned\end{equation}
Applying $$\sigma(R)\rightarrow1,$$
$$D\rightarrow p(n-1)K^2$$
as $R\rightarrow\infty$ in \eqref{In-20}, we have
\begin{equation}\label{In-21}\aligned
|\nabla h|^p(x)=&\varphi(x)|\nabla h|^p(x)
\leq \tilde{G}(x_0)\\
\leq&
\frac{1}{2\Big[\frac{p}{n-1}-\varepsilon_2-\frac{p-2}{p}\varepsilon_3\Big]}\Bigg\{-\Big[\frac{2p}{n-1}(p-1)^{p-1}\lambda-\frac{2}{p}
\frac{[p(n-1)]^{\frac{p}{2}}K^p}{\varepsilon_3^{\frac{p-2}{2}}}\Big]\\
&+\Bigg[\Big(\frac{2p}{n-1}(p-1)^{p-1}\lambda-\frac{2}{p}
\frac{[p(n-1)]^{\frac{p}{2}}K^p}{\varepsilon_3^{\frac{p-2}{2}}}\Big)^2\\
&-4\Big(\frac{p}{n-1}
-\varepsilon_2-\frac{p-2}{p}\varepsilon_3\Big)\Big(\frac{p(p-1)^{2(p-1)}}{n-1}-\varepsilon_1\Big)\lambda^2\Bigg]^{\frac{1}{2}}\Bigg\}.
\endaligned\end{equation}
Taking $\varepsilon_1\rightarrow0$ and $\varepsilon_2\rightarrow0$
in \eqref{In-21}, we have
\begin{equation}\label{HIn-21}\aligned
|\nabla h|^p(x)=&\varphi(x)|\nabla h|^p(x)
\leq \tilde{G}(x_0)\\
\leq&
\frac{1}{2\Big[\frac{p}{n-1}-\frac{p-2}{p}\varepsilon_3\Big]}\Bigg\{-\Big[\frac{2p}{n-1}(p-1)^{p-1}\lambda-\frac{2}{p}
\frac{[p(n-1)]^{\frac{p}{2}}K^p}{\varepsilon_3^{\frac{p-2}{2}}}\Big]\\
&+\Bigg[\Big(\frac{2p}{n-1}(p-1)^{p-1}\lambda-\frac{2}{p}
\frac{[p(n-1)]^{\frac{p}{2}}K^p}{\varepsilon_3^{\frac{p-2}{2}}}\Big)^2\\
&-4\Big(\frac{p}{n-1}
-\frac{p-2}{p}\varepsilon_3\Big)\frac{p(p-1)^{2(p-1)}}{n-1}\lambda^2\Bigg]^{\frac{1}{2}}\Bigg\},
\endaligned\end{equation}
with
\begin{equation}\label{HHIn-21}0<\varepsilon_3<\frac{p^2}{(n-1)(p-2)}.\end{equation}

In particular, letting $R\rightarrow\infty$ and
$\varepsilon_1\rightarrow0$, $\varepsilon_2\rightarrow0$ in
\eqref{HHIn-20}, we obtain
\begin{equation}\label{In-23}\aligned\Big(\frac{2p}{n-1}(p-1)^{p-1}&\lambda-\frac{2}{p}
\frac{[p(n-1)]^{\frac{p}{2}}K^p}{\varepsilon_3^{\frac{p-2}{2}}}\Big)^2\\
&-4\Big(\frac{p}{n-1}
-\frac{p-2}{p}\varepsilon_3\Big)\frac{p(p-1)^{2(p-1)}}{n-1}\lambda^2\geq0.
\endaligned\end{equation}  This shows that
\begin{equation}\label{In-24}
\Big(A^2-4C\Big)\lambda^2-2AB\lambda+B^2\geq0,
\end{equation} where
$$A=\frac{2p}{n-1}(p-1)^{p-1},$$
$$B=\frac{2}{p}
\frac{[p(n-1)]^{\frac{p}{2}}K^p}{\varepsilon_3^{\frac{p-2}{2}}},$$
$$C=\Big(\frac{p}{n-1}
-\frac{p-2}{p}\varepsilon_3\Big)\frac{p(p-1)^{2(p-1)}}{n-1}.$$ When
$p>2$, the coefficient of $\lambda^2$ in \eqref{In-24} positive,
that is,
$$A^2-4C>0.$$ Hence, we have
\begin{equation}\label{In-25}\lambda\leq\frac{2AB-\sqrt{4A^2B^2-4B^2(A^2-4C)}}{2(A^2-4C)}=\frac{B}{A+2\sqrt{C}}.\end{equation}
Minimizing the right hand side of \eqref{In-25} by taking
\begin{equation}\label{In-26}\varepsilon_3=\frac{p^3}{(n-1)(p-1)^2},\end{equation} we obtain
\begin{equation}\label{In-27}\lambda\leq\Big(\frac{(n-1)K}{p}\Big)^p.\end{equation}

In particular, when $\varepsilon_3$ satisfies \eqref{In-26},  the
estimate \eqref{HIn-21} becomes
\begin{equation}\label{In-28}\aligned
|\nabla h|^p(x)\leq&
\frac{(n-1)(p-1)^2}{p}\Bigg\{-\frac{p}{n-1}(p-1)^{p-1}\lambda+(n-1)^{p-1}K^p\Big(\frac{p-1}{p}\Big)^{p-2}\\
&+\Bigg[\frac{p^3(p-1)^{2(p-1)}}{(n-1)^2(p-2)}\lambda^2-\frac{2K^p(p-1)^{2p-3}(n-1)^{p-2}}{p^{p-3}}\lambda\\
& +K^{2p}(n-1)^{2(p-1)}\Big(\frac{p-1}{p}\Big)^{2(p-2)}
\Bigg]^{\frac{1}{2}}\Bigg\}.
\endaligned\end{equation}

This concludes the proof  of Theorem \ref{ThmInt1}.

\section{The Proof of the Case $1<p<2$}

For $1<p<2$, we have
\begin{equation}\label{2-In-3}\min\Big\{2(p-1),\frac{n(p-1)^2}{n-1}\Big\}=\frac{n(p-1)^2}{n-1}.\end{equation}
Then \eqref{In-5} becomes
\begin{equation}\label{2-In-5}\aligned
\mathcal{L}(G)\geq&\frac{p}{n-1}G^2+\frac{2p}{n-1}(p-1)^{p-1}\lambda
G+\frac{p(p-1)^{2(p-1)}}{n-1}\lambda^2\\
&-p(n-1)K^2G^{\frac{2(p-1)}{p}}+\frac{n(p-1)^2}{(n-1)p}\frac{|\nabla
G|^2}{G^{\frac{2}{p}}}
\\
&+\Bigg[\frac{2(p-1)}{n-1}\Big(1+\lambda
G^{-1}(p-1)^{p-1}\Big)-p\Bigg]\langle\nabla h,\nabla G\rangle
|\nabla h|^{p-2}.
\endaligned\end{equation}
Similarly, we assume that $\tilde{G}=\varphi G$ achieve its maximal
value at $x_0$, where $\varphi$ is defined in Section 2. Then at the
point $x_0$, we also have that \eqref{In-7} is valid. By a direct
calculation and noticing $1<p<2$, we can deduce from \eqref{2-In-5}
that
\begin{equation}\label{2-In-8}\aligned
\varphi\mathcal{L} G \geq&\frac{p}{n-1}\
\frac{\tilde{G}^2}{\varphi}+\frac{2p}{n-1}(p-1)^{p-1}\lambda
\tilde{G}+\frac{p(p-1)^{2(p-1)}}{n-1}\lambda^2\varphi\\
&-p(n-1)K^2\frac{\tilde{G}^{\frac{2(p-1)}{p}}}{\varphi^{\frac{p-2}{p}}}
+\frac{n(p-1)^2}{(n-1)p}\frac{|\nabla\varphi|^2}{\varphi}\Big(\frac{\tilde{G}}{\varphi}\Big)^{\frac{2(p-1)}{p}}
\\
&+\frac{\tilde{G}}{\varphi}\Bigg(p-\frac{2(p-1)}{n-1}
\Bigg)\langle\nabla h,\nabla \varphi\rangle
|\nabla h|^{p-2}\\
&-\frac{2}{n-1}(p-1)^p\lambda \langle\nabla h,\nabla \varphi\rangle
|\nabla h|^{p-2}\\
\geq&\frac{p}{n-1}\
\frac{\tilde{G}^2}{\varphi}+\frac{2p}{n-1}(p-1)^{p-1}\lambda
\tilde{G}+\frac{p(p-1)^{2(p-1)}}{n-1}\lambda^2\varphi\\
&-p(n-1)K^2\frac{\tilde{G}^{\frac{2(p-1)}{p}}}{\varphi^{\frac{p-2}{p}}}
+\frac{n(p-1)^2}{(n-1)p}\frac{|\nabla\varphi|^2}{\varphi}\Big(\frac{\tilde{G}}{\varphi}\Big)^{\frac{2(p-1)}{p}}
\\
&-\Big(p-\frac{2(p-1)}{n-1}
\Big)|\nabla\varphi|\Big(\frac{\tilde{G}}{\varphi}\Big)^{\frac{2p-1}{p}}-\frac{2}{n-1}(p-1)^p\lambda
|\nabla\varphi|\Big(\frac{\tilde{G}}{\varphi}\Big)^{\frac{p-1}{p}}.
\endaligned\end{equation}
On the other hand,
\begin{equation}\label{2-In-9}\aligned
2\langle\nabla\varphi,A(\nabla G)\rangle|\nabla
h|^{p-2}=&2\Bigg[\langle\nabla\varphi, \nabla G
\rangle+(p-2)\frac{\langle\nabla \varphi,\nabla
h\rangle\langle\nabla h,\nabla G\rangle}{|\nabla
h|^{2}}\Bigg]|\nabla h|^{p-2}\\
=&-2G\Bigg[\frac{|\nabla\varphi|^2}{\varphi}+(p-2)\frac{\langle\nabla
h,\nabla \varphi\rangle^2}{\varphi|\nabla h|^{2}}\Bigg]|\nabla
h|^{p-2}\\
\geq&-2\frac{|\nabla\varphi|^2}{\varphi}\Big(\frac{\tilde{G}}{\varphi}\Big)^{\frac{2(p-1)}{p}}.
\endaligned\end{equation}
For $1<p<2$, we also have
\begin{equation}\label{2-In-10}\aligned
\mathcal{L}\varphi=&{\rm div}(|\nabla
h|^{p-2}A(\nabla\varphi))\\
=&\varphi_{ij}A^{ij}\Big(\frac{\tilde{G}}{\varphi}\Big)^{\frac{p-2}{p}}-(p-2)\frac{\langle\nabla\varphi,\nabla
h\rangle}{|\nabla
h|^2}\Big(G+(p-1)^{p-1}\lambda\Big)\\
&+2\Big(\frac{2}{p}-1\Big)\frac{|\nabla\varphi|^2}{\varphi}\Big(\frac{\tilde{G}}{\varphi}\Big)^{\frac{p-2}{p}}
+\frac{2(p-2)}{p}\frac{\langle\nabla\varphi,\nabla
h\rangle^2}{|\nabla
h|^2\varphi}\Big(\frac{\tilde{G}}{\varphi}\Big)^{\frac{p-2}{p}}.
\endaligned\end{equation}
Therefore, we obtain from \eqref{2-In-10}
\begin{equation}\label{2-In-11}\aligned
G\mathcal{L}\varphi=&\varphi_{ij}A^{ij}\Big(\frac{\tilde{G}}{\varphi}\Big)^{\frac{2(p-1)}{p}}+(2-p)\langle\nabla\varphi,\nabla
h\rangle\Big(\frac{\tilde{G}}{\varphi}\Big)^{\frac{2(p-1)}{p}}\\
&+(2-p)(p-1)^{p-1}\lambda\langle\nabla\varphi,\nabla
h\rangle\Big(\frac{\tilde{G}}{\varphi}\Big)^{\frac{p-2}{p}}\\
&+2\Big(\frac{2}{p}-1\Big)\frac{|\nabla\varphi|^2}{\varphi}\Big(\frac{\tilde{G}}{\varphi}\Big)^{\frac{2(p-1)}{p}}
-\frac{2(2-p)}{p}\frac{\langle\nabla\varphi,\nabla
h\rangle^2}{|\nabla
h|^2\varphi}\Big(\frac{\tilde{G}}{\varphi}\Big)^{\frac{2(p-1)}{p}}\\
\geq&\varphi_{ij}A^{ij}\Big(\frac{\tilde{G}}{\varphi}\Big)^{\frac{2(p-1)}{p}}
-(2-p)|\nabla\varphi|\Big(\frac{\tilde{G}}{\varphi}\Big)^{\frac{2p-1}{p}}\\
&-(2-p)(p-1)^{p-1}\lambda|\nabla\varphi|\Big(\frac{\tilde{G}}{\varphi}\Big)^{\frac{p-1}{p}}
+2\Big(\frac{2}{p}-1\Big)\frac{|\nabla\varphi|^2}{\varphi}\Big(\frac{\tilde{G}}{\varphi}\Big)^{\frac{2(p-1)}{p}}\\
&-\frac{2(2-p)}{p}\frac{\langle\nabla\varphi,\nabla
h\rangle^2}{|\nabla
h|^2\varphi}\Big(\frac{\tilde{G}}{\varphi}\Big)^{\frac{2(p-1)}{p}}\\
\geq&\varphi_{ij}A^{ij}\Big(\frac{\tilde{G}}{\varphi}\Big)^{\frac{2(p-1)}{p}}
-(2-p)|\nabla\varphi|\Big(\frac{\tilde{G}}{\varphi}\Big)^{\frac{2p-1}{p}}\\
&-(2-p)(p-1)^{p-1}\lambda|\nabla\varphi|\Big(\frac{\tilde{G}}{\varphi}\Big)^{\frac{p-1}{p}}.
\endaligned\end{equation}
With the help of \eqref{2-In-8}, \eqref{2-In-9} and \eqref{2-In-11},
we obtain
\begin{equation}\label{2-In-12}\aligned
\mathcal{L}(\tilde{G})=&\varphi\mathcal{L}
G+G\mathcal{L}\varphi+2\langle\nabla\varphi,A(\nabla
G)\rangle|\nabla
h|^{p-2}\\
\geq&\frac{p}{n-1}\
\frac{\tilde{G}^2}{\varphi}+\frac{2p}{n-1}(p-1)^{p-1}\lambda
\tilde{G}+\frac{p(p-1)^{2(p-1)}}{n-1}\lambda^2\varphi\\
&+\Bigg[-p(n-1)K^2\varphi
+\frac{n(p-1)^2}{(n-1)p}\frac{|\nabla\varphi|^2}{\varphi}-2\frac{|\nabla\varphi|^2}{\varphi}\\
&+\varphi_{ij}A^{ij}\Bigg]\Big(\frac{\tilde{G}}{\varphi}\Big)^{\frac{2(p-1)}{p}}
-\frac{2(n-p)}{n-1}|\nabla\varphi|\Big(\frac{\tilde{G}}{\varphi}\Big)^{\frac{2p-1}{p}}\\
&-\frac{(2-p)n+3p-4}{n-1}(p-1)^{p-1}\lambda|\nabla\varphi|\Big(\frac{\tilde{G}}{\varphi}\Big)^{\frac{p-1}{p}}.
\endaligned\end{equation}
Using \eqref{H-In-6} and the Cauchy inequality, we have
\begin{equation}\label{2-In-13}\aligned
-\frac{(2-p)n+3p-4}{n-1}&(p-1)^{p-1}\lambda|\nabla\varphi|\Big(\frac{\tilde{G}}{\varphi}\Big)^{\frac{p-1}{p}}\\
\geq&-\sqrt{10}\frac{(2-p)n+3p-4}{n-1}(p-1)^{p-1}\frac{\lambda}{R}\Big(\frac{\tilde{G}}{\varphi}\Big)^{\frac{p-1}{p}}\varphi^{\frac{1}{2}}\\
\geq&-\varepsilon_1\lambda^2
-\frac{5(p-1)^{2(p-1)}\frac{[(2-p)n+3p-4]^2}{(n-1)^2}}{2\varepsilon_1R^2}
\Big(\frac{\tilde{G}}{\varphi}\Big)^{\frac{2(p-1)}{p}}
\endaligned\end{equation} and
\begin{equation}\label{2-In-14}\aligned
-\frac{2(n-p)}{n-1}|\nabla\varphi|\Big(\frac{\tilde{G}}{\varphi}\Big)^{\frac{2p-1}{p}}
=&-\frac{2(n-p)}{n-1}|\nabla\varphi|\Big(\frac{\tilde{G}}{\varphi}\Big)^{\frac{p-1}{p}}\Big(\frac{\tilde{G}}{\varphi}\Big)\\
\geq&-\varepsilon_2\Big(\frac{\tilde{G}}{\varphi}\Big)^2-\frac{4(n-p)^2}{(n-1)^2}
\frac{5\varphi}{2\varepsilon_2R^2}\Big(\frac{\tilde{G}}{\varphi}\Big)^{\frac{2(p-1)}{p}},
\endaligned\end{equation}
where $\varepsilon_1, \varepsilon_2$ are two positive constants.
Inserting \eqref{2-In-13} and \eqref{2-In-14} into \eqref{2-In-12},
we derive
\begin{equation}\label{2-In-15}\aligned
\mathcal{L}(\tilde{G})\geq&\Big[\frac{p}{n-1}\varphi-\varepsilon_2\Big]\Big(\frac{\tilde{G}}{\varphi}\Big)^2+\frac{2p}{n-1}(p-1)^{p-1}\lambda\varphi
\frac{\tilde{G}}{\varphi}+\Big[\frac{p(p-1)^{2(p-1)}}{n-1}\varphi-\varepsilon_1\Big]\lambda^2\\
&-\Bigg\{p(n-1)K^2\varphi
+\frac{-np^2+2(2n-1)p-n}{p(n-1)}\frac{|\nabla\varphi|^2}{\varphi}\\
&-\varphi_{ij}A^{ij}+\frac{5(p-1)^{2(p-1)}\frac{[(2-p)n+3p-4]^2}{(n-1)^2}\varphi}{2\varepsilon_1R^2}
+\frac{4(n-p)^2}{(n-1)^2}\frac{5\varphi}{2\varepsilon_2R^2}\Bigg\}\Big(\frac{\tilde{G}}{\varphi}\Big)^{\frac{2(p-1)}{p}}\\
\geq&\Big[\frac{p}{n-1}\varphi-\varepsilon_2\Big]\Big(\frac{\tilde{G}}{\varphi}\Big)^2+\frac{2p}{n-1}(p-1)^{p-1}\lambda\varphi
\frac{\tilde{G}}{\varphi}+\Big[\frac{p(p-1)^{2(p-1)}}{n-1}\varphi-\varepsilon_1\Big]\lambda^2\\
&-\overline{D}\Big(\frac{\tilde{G}}{\varphi}\Big)^{\frac{2(p-1)}{p}},
\endaligned\end{equation}
where $-np^2+2(2n-1)p-n>0$ and
$$\aligned \overline{D}=&p(n-1)K^2\varphi
+\frac{-np^2+2(2n-1)p-n}{p(n-1)}\frac{10}{R^2}\\
&+80(n+p-2)\frac{1+KR}{R^2}+\frac{40}{R^2}+\frac{5(p-1)^{2(p-1)}\frac{[(2-p)n+3p-4]^2}{(n-1)^2}\varphi}{2\varepsilon_1R^2}\\
&+\frac{4(n-p)^2}{(n-1)^2}\frac{5\varphi}{2\varepsilon_2R^2}.\endaligned$$
By virtue of the Holder inequality, one has
\begin{equation}\label{2-In-16}\aligned
\overline{D}\Big(\frac{\tilde{G}}{\varphi}\Big)^{\frac{2(p-1)}{p}}
\leq&\frac{2-p}{p}\varepsilon_3^{\frac{2(p-1)}{2-p}}+\frac{2(p-1)}{p}
\frac{\overline{D}^{\frac{p}{2(p-1)}}}{\varepsilon_3}\frac{\tilde{G}}{\varphi}.
\endaligned\end{equation}
Thus, from \eqref{2-In-15} we deduces to
\begin{equation}\label{2-In-17}\aligned
\mathcal{L}(\tilde{G})\geq&
\Big[\frac{p}{n-1}\varphi-\varepsilon_2\Big]\Big(\frac{\tilde{G}}{\varphi}\Big)^2\\
&+\Big[\frac{2p}{n-1}(p-1)^{p-1}\lambda\varphi-\frac{2(p-1)}{p}
\frac{\overline{D}^{\frac{p}{2(p-1)}}}{\varepsilon_3}\Big]
\frac{\tilde{G}}{\varphi}\\
&+\Big[\frac{p(p-1)^{2(p-1)}}{n-1}\varphi-\varepsilon_1\Big]\lambda^2
-\frac{2-p}{p}\varepsilon_3^{\frac{2(p-1)}{2-p}},
\endaligned\end{equation}
which gives that at the point $x_0$,
\begin{equation}\label{2-In-18}\aligned
0\geq&
\Big[\frac{p}{n-1}\varphi-\varepsilon_2\Big]\tilde{G}^2\\
&+\Big[\frac{2p}{n-1}(p-1)^{p-1}\lambda\varphi^2-\frac{2(p-1)}{p}
\frac{\overline{D}^{\frac{p}{2(p-1)}}}{\varepsilon_3}\varphi\Big]
\tilde{G}\\
&+\Big[\frac{p(p-1)^{2(p-1)}}{n-1}\varphi^3-\varepsilon_1\varphi^2\Big]\lambda^2-\frac{2-p}{p}\varepsilon_3^{\frac{2(p-1)}{2-p}}\varphi^2.
\endaligned\end{equation}

Similarly, we define
$$\sigma(R):=\frac{\sup\limits_{B_p(\frac{R}{2})}|\nabla h|^p(x)}{\sup\limits_{B_p(R)}|\nabla
h|^p(x)}.$$ Then, it is easy to see that
$$\sigma(R)\leq\varphi(x_0)\leq1,$$
and at the point $x_0$, \eqref{2-In-18} gives
\begin{equation}\label{2-In-19}\aligned
0\geq&
\Big[\frac{p}{n-1}\sigma(R)-\varepsilon_2\Big]\tilde{G}^2\\
&+\Big[\frac{2p}{n-1}(p-1)^{p-1}\lambda\sigma(R)^2-\frac{2(p-1)}{p}
\frac{\overline{D}^{\frac{p}{2(p-1)}}}{\varepsilon_3}\Big]
\tilde{G}\\
&+\Big[\frac{p(p-1)^{2(p-1)}}{n-1}\sigma(R)^3-\varepsilon_1\Big]\lambda^2
-\frac{2-p}{p}\varepsilon_3^{\frac{2(p-1)}{2-p}},
\endaligned\end{equation}
 as long as
$$\frac{p}{n-1}\sigma(R)-\varepsilon_2>0.$$
By using the inequality
$$ax^2+bx+c\leq0$$ with $a>0$, then
$$x\leq\frac{-b+\sqrt{b^2-4ac}}{2a},$$
we obtain from \eqref{2-In-19}
\begin{equation}\label{H2-In-20}\aligned
&\Big[\frac{2p}{n-1}(p-1)^{p-1}\lambda\sigma(R)^2-\frac{2(p-1)}{p}
\frac{\overline{D}^{\frac{p}{2(p-1)}}}{\varepsilon_3}\Big]^2\\
&-4\Big[\frac{p}{n-1}\sigma(R)-\varepsilon_2\Big]\Big\{\Big[\frac{p(p-1)^{2(p-1)}}{n-1}\sigma(R)^3-\varepsilon_1\Big]\lambda^2
-\frac{2-p}{p}\varepsilon_3^{\frac{2(p-1)}{2-p}}\Big\}\\
&\geq0
\endaligned\end{equation}
and
\begin{equation}\label{2-In-20}\aligned
\tilde{G}(x_0)\leq&
\frac{1}{2\Big[\frac{p}{n-1}\sigma(R)-\varepsilon_2\Big]}
\Bigg\{-\Big[\frac{2p}{n-1}(p-1)^{p-1}\lambda\sigma(R)^2-\frac{2(p-1)}{p}
\frac{\overline{D}^{\frac{p}{2(p-1)}}}{\varepsilon_3}\Big]\\
&+\Bigg\{\Big(\frac{2p}{n-1}(p-1)^{p-1}\lambda\sigma(R)^2-\frac{2(p-1)}{p}
\frac{\overline{D}^{\frac{p}{2(p-1)}}}{\varepsilon_3}\Big)^2\\
&-4\Big(\frac{p}{n-1}\sigma(R)
-\varepsilon_2\Big)\Bigg[\Big(\frac{p(p-1)^{2(p-1)}}{n-1}\sigma(R)^3-\varepsilon_1\Big)\lambda^2
-\frac{2-p}{p}\varepsilon_3^{\frac{2(p-1)}{2-p}}\Bigg]\Bigg\}^{\frac{1}{2}}\Bigg\}
\endaligned\end{equation}
Noting $$\sigma(R)\rightarrow1,$$
$$\overline{D}\rightarrow p(n-1)K^2$$
as $R\rightarrow\infty$, hence \eqref{2-In-20} becomes
\begin{equation}\label{2-In-21}\aligned
|\nabla h|^p(x)\leq&\frac{n-1}{2p}
\Bigg\{-\Big[\frac{2p}{n-1}(p-1)^{p-1}\lambda-\frac{2(p-1)}{p}
\frac{[p(n-1)K^2]^{\frac{p}{2(p-1)}}}{\varepsilon_3}\Big]\\
&+\Bigg[\Big(\frac{2p}{n-1}(p-1)^{p-1}\lambda-\frac{2(p-1)}{p}
\frac{[p(n-1)K^2]^{\frac{p}{2(p-1)}}}{\varepsilon_3}\Big)^2\\
&-\frac{4p}{n-1}\Bigg(\frac{p(p-1)^{2(p-1)}}{n-1}\lambda^2
-\frac{2-p}{p}\varepsilon_3^{\frac{2(p-1)}{2-p}}\Bigg)\Bigg]^{\frac{1}{2}}\Bigg\}
\endaligned\end{equation} by letting $\varepsilon_1\rightarrow0$ and
$\varepsilon_2\rightarrow0$.

In particular, letting $R\rightarrow\infty$ and
$\varepsilon_1,\varepsilon_2\rightarrow0$ in \eqref{H2-In-20}, we
obtain
\begin{equation}\label{2-In-22}\aligned&\Big(\frac{2p}{n-1}(p-1)^{p-1}\lambda-\frac{2(p-1)}{p}
\frac{[p(n-1)K^2]^{\frac{p}{2(p-1)}}}{\varepsilon_3}\Big)^2\\
&-\frac{4p}{n-1}\Bigg(\frac{p(p-1)^{2(p-1)}}{n-1}\lambda^2
-\frac{2-p}{p}\varepsilon_3^{\frac{2(p-1)}{2-p}}\Bigg)\geq0,
\endaligned\end{equation}
which shows that
\begin{equation}\label{2-In-24}
\Big(\tilde{A}^2-4\tilde{C}\Big)\lambda^2-2\tilde{A}\tilde{B}\lambda+\tilde{B}^2+4\tilde{D}\geq0,
\end{equation} where
$$\tilde{A}=\frac{2p}{n-1}(p-1)^{p-1},$$
$$\tilde{B}=\frac{2(p-1)}{p}
\frac{[p(n-1)K^2]^{\frac{p}{2(p-1)}}}{\varepsilon_3},$$
$$\tilde{C}=\frac{p^2(p-1)^{2(p-1)}}{(n-1)^2},$$
$$\tilde{D}=\frac{2-p}{n-1}\varepsilon_3^{\frac{2(p-1)}{2-p}}.$$
Since the coefficient of $\lambda^2$ is zero, that is,
$$\tilde{A}^2-4\tilde{C}=0.$$
Hence, we have
\begin{equation}\label{2-In-25}
\lambda\leq\frac{\tilde{B}^2+4\tilde{D}}{2\tilde{A}\tilde{B}}.
\end{equation}
Minimizing the right hand side of \eqref{2-In-25} by taking
\begin{equation}\label{2-In-26}\varepsilon_3=\Big(\frac{(n-1)(p-1)^2}{p^3}\Big)^{\frac{2-p}{2}}
[p(n-1)K^2]^{\frac{p(2-p)}{2(p-1)}} ,\end{equation} with the
relationship
$$\frac{4p}{n-1}\varepsilon_3^{\frac{2(p-1)}{2-p}}=\tilde{B}^2,$$
we obtain
\begin{equation}\label{2-In-27}\lambda\leq\Big(\frac{(n-1)K}{p}\Big)^p.\end{equation}

In particular, when $\varepsilon_3$ satisfies \eqref{2-In-26},  the
estimate \eqref{2-In-21} becomes
\begin{equation}\label{2-In-28}\aligned
|\nabla h|^p(x)\leq&\frac{n-1}{p}
\Bigg\{-\frac{p}{n-1}(p-1)^{p-1}\lambda+(n-1)^{p-1}p^{2-p}(p-1)^{p-1}K^p
\\
&+\Bigg[-\frac{2(p-1)^{2(p-1)}p^{3-p}K^p}{(n-1)^{2-p}}\lambda+2(n-1)^{2(p-1)}(p-1)^{2(p-1)}p^{3-2p}K^{2p}\Bigg]^{\frac{1}{2}}\Bigg\}
\endaligned\end{equation}

This concludes the proof  of Theorem \ref{ThmInt2}.

\bibliographystyle{Plain}

\end{document}